\PassOptionsToPackage{square}{natbib}  
\documentclass[final,3p,times,authoryear]{elsarticle}
\usepackage{amsmath}
\usepackage{amssymb}
\usepackage{amsmath,amssymb}
\usepackage[authoryear]{natbib}
\usepackage[utf8]{inputenc}
\usepackage[english]{babel}
\usepackage{mathtools}
\usepackage{xcolor}
\usepackage[overload]{empheq}
\usepackage{ctable}
\usepackage{cases,eufrak,mathrsfs} 
\definecolor{mygreen}{RGB}{10,150,10}

\newcommand\myblue[1]{\textcolor{black}{#1}}

\usepackage{algorithm}
\usepackage{algpseudocode}
\usepackage{multicol}
\usepackage{multirow}
\usepackage{tabu,threeparttable}
\usepackage{mathrsfs}
\journal{Intl. of  Journal of Computer Mathematics}

\begin{document}
	{

%
%
%
\begin{frontmatter}

\title{Preconditioned Linear Solves for Parametric Model Order Reduction}  

\author{Navneet Pratap Singh 
 and Kapil Ahuja}
\address{Computer Science and Engineering, Indian Institute of Technology Indore, India}




\begin{abstract}
There exist many classes of algorithms for computing reduced-order models of parametric dynamical systems, commonly termed as  parametric model order reduction algorithms. The main computational cost of these algorithms is in solving sequences of very large and sparse linear systems of equations, which are predominantly dependent on slowly varying parameter values. We focus on efficiently solving these linear systems, arising while reducing second-order linear dynamical systems, by iterative methods with appropriate preconditioners. We propose that the choice of underlying iterative solver is problem dependent. Since for many parametric model order reduction algorithms, the linear systems right-hand-sides are available together, we propose the use of block variant of the underlying iterative method.

Due to constant increase in the input model size and the number of parameters in it, computing a preconditioner in a parallel setting is increasingly becoming a norm. Since, Sparse Approximate Inverse (SPAI) preconditioner is a general preconditioner that can be naturally parallelized, we propose its use. Our most novel contribution is a technique to cheaply update the SPAI preconditioner, while solving the parametrically changing linear systems.  We support our proposed theory by numerical experiments where we first show that using a block variant of the underlying iterative solver saves 80\% of the computation time over the non-block version. Further, and more importantly, SPAI with updates saves 70\% of the time over SPAI without updates.
\end{abstract}

\begin{keyword}
	Parametric Model Order Reduction, Parametrically Dependent Linear Systems, Iterative Methods, SPAI Preconditioner, and Preconditioner Updates.
\MSC[2010] 34C20 \sep 65F10
\end{keyword}

\end{frontmatter}


\section{Introduction}
\label{sec:Intro}
  
Dynamical systems arise while modelling of many engineering and scientific applications \cite{Ogata2001, ANTOULAS200419}. These dynamical systems depend upon parameters, which vary with different design stages or computer experiments. Substantial work has been done for the first  order linear systems \cite{Ogata2001, ANTOULAS200419}, and hence, we focus on  second order linear systems here. Higher order systems can also be looked at, which is part of our future work.

 A parameterized second order linear dynamical system is usually of the form 
\begin{align}\label{eq:Linear_param_system}
\begin{split}
M(p_j)\ddot{x}(t) + D(p_j)\dot{x}(t) + K(p_j)x(t)  = Bu(t), \\ 
y(t)  = C_1(p_j) \dot{x}(t) + C_2(p_j) x(t),
\end{split}
\end{align}
where $M(p_j), D(p_j), K(p_j) \in \mathcal{R}^{n \times n}$, $B \in \mathcal{R}^{n \times d_I}$, $C_1, C_2 \in \mathcal{R}^{d_O \times n}$ and $p_j$ for $j = 1, \ldots w$, are the parameters.
Also, $x(t) \colon \mathcal{R} \rightarrow \mathcal{R}^n$ is the vector of all states,  $u(t) \colon \mathcal{R} \rightarrow \mathcal{R}^{d_I}$ and $y(t) \colon \mathcal{R} \rightarrow \mathcal{R}^{d_O}$ are the inputs and the outputs of the system, respectively. 

Parameterized dynamical systems of this type are usually  very large in size. Solving such systems by traditional simulation methods is often very time-consuming. In these systems,  the parameters are to be provided as fixed values that cannot be changed during simulation. Moreover, since  many runs with differing parameters values are required \cite{lihong2014}, the simulation process is enormous.

Model order reduction, traditionally developed for non-parametric systems \cite{Meyer1996, Beattie2005, BONIN20161}, is a very popular technique to overcome such issues. A reduced system can be derived by model order reduction that can then be used for simulation instead of the full system. This process often saves substantial simulation time. Model order reduction for parametric systems, parametric model order reduction, preserves the parameters of the original system as a symbolic quantities in the reduced system. Whenever there is a change in the parameters, we need not recompute the new reduced system. Instead, we simply use the changed parameters while solving the reduced system. 



Many algorithms exist for model order reduction of parametrized second order linear dynamical systems. Some common ones are as follows: 
	 Robust Algorithm for Parametric Model Order Reduction (RPMOR) \cite{lihong2014}, which is based on moment matching; Data-Driven Parametrized Model Reduction algorithm in the Loewner Framework (PMOR-L) \cite{Ionita2014}; and Parametric Balanced Truncation Model Reduction algorithm (PBTMR) \cite{sonN2017}, which is based on a Greedy approach. 
Next, we summarize these three algorithms, and then abstract out the computational bottleneck step of solving linear systems.  

\subsection{Robust Algorithm for Parametric Model Order Reduction}
RPMOR \cite{lihong2014} is a projection based model order reduction algorithm and is mainly used for the reduction of parametric first and second order linear dynamical systems. 
Here, the state variable $x(t)$ is projected onto a smaller dimensional subspace. Let  $V \in \mathcal{R}^{n \times r}$ be a projection matrix determined by RPMOR. Using $x(t) \approx V \hat{x}(t)$, $C_1(p_j) = 0$, and $C_2(p_j) = C$ in (\ref{eq:Linear_param_system}), we obtain the following system: 
\begin{align*}
\begin{split}
M(p_j)V\ddot{\hat{x}}(t) + D(p_j)V\dot{\hat{x}}(t) + K(p_j)V\hat{x}(t) - Bu(t) & = r(t), \\
\hat{y}(t) & = C^TV \hat{x}(t),
\end{split}
\end{align*}
where $r(t)$ is the residual after projection. Applying the Galerkin approach by multiplying $V^T$ in the first equation above we get 
\begin{align*}
\begin{split}
V^T\left(M(p_j)V\ddot{\hat{x}}(t) + D(p_j)V\dot{\hat{x}}(t) + K(p_j)V\hat{x}(t) - Bu(t)\right) & = 0, \\
\hat{y}(t)  & = C^TV \hat{x}(t) \\
or \qquad \qquad \qquad \qquad \qquad
\end{split}
\end{align*}
\begin{align}\label{eq:red_sys}
\begin{split}
\hat{M}(p_j)\ddot{\hat{x}}(t) + \hat{D}(p_j)\dot{\hat{x}}(t) + \hat{K}(p_j)\hat{x}(t)  - \hat{B}u(t) & = 0, \\
\hat{y}(t) & = \hat{C}^T\hat{x}(t),
\end{split}
\end{align}
where $\hat{M}(p_j), \ \hat{D}(p_j),  \ \hat{K}(p_j) \ \in \mathcal{R}^{r \times r}, \hat{B} \ \in \mathcal{R}^{r \times d_I},  \hat{C} \ \in \mathcal{R}^{d_O \times r}$, and $r << n$. We want  $\hat{y}(t)$ should be nearly equal to $y(t)$ for all acceptable inputs.

The projection matrix $V$ can be determined by many ways. One common way is by moment matching \cite{Benner2015,lihong2014,Serkan2002}, which is discussed next.  
In the frequency domain, (\ref{eq:Linear_param_system}) with $C_1(p_j) = 0$ and $C_2(p_j) = C$ is given as
\begin{align*}
\begin{split}
\left(s^2 M(p_j) + s D(p_j) + K(p_j)\right)x(s) & = B u(s), \\ y(s) & = C^Tx(s),
\end{split}
\end{align*}
where $s$ is the new parameter (frequency parameter corresponding to time $t$). 
The above equation can also be rewritten as
\begin{align}\label{eq:gen_par_freq}
\begin{split}
A\left(s, p_j\right)x(s) & = Bu(s), \\ y(s) & = C^Tx(s),
\end{split}
\end{align}
where  $A\left(s, p_j\right) \in \mathcal{R}^{n \times n}$ is the  parametrized matrix.
Next, the system in (\ref{eq:gen_par_freq}) is transformed to an affine form as
\begin{align}\label{eq:affine_system}
\begin{split}
\left(A_0 + \tilde{s}_1 A_1 + \cdots + \tilde{s}_w A_w + \tilde{s}_{w+1} A_{w+1}\right)x(s)  & = Bu(s), \\ y(s) & = C^Tx(s),
\end{split}
\end{align}
where $A_0, A_1, \ldots, A_{w+1} \in \mathcal{R}^{n \times n}$, the new parameters $\tilde{s}_1$ and $\tilde{s}_j \ (\text{for} \  j = 2,\ldots, w+1)$ are some functions (polynomial, rational, etc.) of the parameters $s$ and $p_j$, respectively. 
Next, 
the state $x(s)$ in (\ref{eq:affine_system}) is computed at initial expansion point (from here onwards we represent set of parameters as an expansion point) $\tilde{\tilde{s}}_1 = \left(\tilde{s}^1_1, \ldots, \tilde{s}^1_w, \tilde{s}^1_{w+1}\right)$ as
\begin{align}\label{eq:momt_exp}
x(s) & = \left[I - \left(\sigma_1 M_1 + \ldots + \sigma_w M_w + \sigma_{w+1} M_{w+1} \right)\right]^{-1}\left(A{(1)}\right)^{-1}Bu(s),
\end{align}
where $\sigma_j = \tilde{s}_j-\tilde{s}_j^1, M_{j} = -\left(A{(1)}\right)^{-1}A_j \ \text{for} \ j = 1, 2, \ldots, w+1,$ and 
\begin{align}\label{eq:param_mat1}
A{(1)} = A_0 + \tilde{s}_1^1  A_1 + \tilde{s}_2^1 A_2 + \ldots + \tilde{s}_{w+1}^1 A_{w+1}.
\end{align}
Applying Taylor series expansion on (\ref{eq:momt_exp}) we get
\begin{align}\label{eq:momt-org}
\begin{split}
x(s)  & = \sum_{h = 0}^{\infty} \left[\sigma_1M_1 + \ldots + \sigma_w M_w + \sigma_{w+1} M_{w+1} \right]^h \tilde{B} u(s), \
\\
& = \sum_{h = 0}^{\infty} x^{(h)}(\sigma_1,\ldots,\sigma_w, \sigma_{w+1}) u(s),
\end{split}
\end{align}
where 
\begin{align*}
\tilde{B} & = \left(A(1)\right)^{-1}B, \\
x^{(0)}(\sigma_1,\ldots,\sigma_w, \sigma_{w+1}) & = \tilde{B}, \\
x^{(1)}(\sigma_1,\ldots,\sigma_w, \sigma_{w+1}) & = \left[\sigma_1M_1 + \ldots + \sigma_wM_w + \sigma_{w+1}M_{w+1}\right] x^{(0)}(\sigma_1,\ldots,\sigma_w, \sigma_{w+1}), \\
& \vdots \\
x^{(h)}(\sigma_1,\ldots,\sigma_w, \sigma_{w+1}) & = \left[\sigma_1M_1 + \ldots + \sigma_wM_w + \sigma_{w+1}M_{w+1}\right] x^{(h-1)}(\sigma_1,\ldots,\sigma_w, \sigma_{w+1}).
\end{align*}
Here, $x^{(h)}(\sigma_1,\ldots,\sigma_w, \sigma_{w+1})$ is called the $h^{th}$-order system moment at $(\sigma_1,\ldots,\sigma_w, \sigma_{w+1})$. Similarly for the reduced system (\ref{eq:red_sys}), the state variable can be written as
\begin{align}\label{eq:momt_red}
\begin{split}
\hat{x}(s) = 
\sum_{h = 0}^{\infty} \hat{x}^{(h)}(\sigma_1,\ldots,\sigma_w, \sigma_{w+1}) u(s).
\end{split}
\end{align}
In the reduced system, the $h^{th}$-order system moment $\hat{x}^{(h)}(\sigma_1,\ldots,\sigma_w, \sigma_{w+1})$ is defined similar to $x^{(h)}(\sigma_1,\ldots,\sigma_w, \sigma_{w+1})$. The goal of moment matching approach is to find a reduced system such that the first few moments of (\ref{eq:momt-org}) and (\ref{eq:momt_red}) are matched. 
This provides the orthogonal projection matrix $V$. 
The columns of $V$ are given by $\it{span} \left \{x^{(0)}(\sigma_1,\ldots,\sigma_w, \sigma_{w+1}), \right.  \\ \left. x^{(1)}(\sigma_1,\ldots,\sigma_w, \sigma_{w+1}), \ldots,  x^{(h)}  (\sigma_1,\ldots,\sigma_w, \sigma_{w+1}) \right \}$, where $h \in \mathcal{R}$. 

After obtaining the first few columns of $V$ corresponding to the current expansion point, the above process is repeated with new set of  expansion point $\tilde{\tilde{s}}_2 = \left(\tilde{s}^2_1, \ldots, \tilde{s}^2_w, \tilde{s}^2_{w+1}\right)$, and we get
\begin{align}\label{eq:param_mat2}
A{(2)} = A_0 + \tilde{s}_1^2  A_1 + \tilde{s}_2^2 A_2 + \ldots + \tilde{s}_{w+1}^2 A_{w+1}.
\end{align}
A similar process can be used for the afterwards set of expansion points $\tilde{\tilde{s}}_{\ell} = \left(\tilde{s}^{\ell}_1, \ldots, \tilde{s}^{\ell}_w, \right. \left. 
\tilde{s}^{\ell}_{w+1}\right)$ for $\ell = 1, 2, \ldots,   \mathfrak{z}$, and we get
\begin{align}\label{eq:param_matl}
A{(\ell)} = A_0 + \tilde{s}_1^{\ell}  A_1 + \tilde{s}_2^{\ell} A_2 + \ldots + \tilde{s}_{w+1}^{\ell} A_{w+1}.
\end{align}

The work in \cite{lihong2014} proposes a popular algorithm based upon this theory (Algorithm 6.1 in \cite{lihong2014}). 

\subsection{Other Parametric Model Order Reduction Algorithms}
    
	
	\myblue{PMOR-L \cite{Ionita2014} is a Loewner framework based  model order reduction algorithm and is used for the model reduction of all types of parametric dynamical systems (linear-nonlinear; first order-higher orders). One important step here is computing a matrix called the Loewner matrix. The computation of this matrix requires computation of the transfer function of the dynamical system. Since we are focussing on second order linear systems, this function for (\ref{eq:Linear_param_system}) with $C_1({p_j}) = 0$ and $C_2(p_j) = C(p_j)$ is given by
	\begin{align} \label{eq:loewner_eq}
	H(s_k,p_j) = C(p_j)^T A(s_k,p_j)^{-1}B(p_j) \ \ \text{for} \ k = 1,\ldots, v, \ \text{and} \  j = 1, \ldots, w,
	\end{align}    
	where 
	\begin{align}\label{eq:loewner_eq2}
	A(s_k,p_j) = \left(s_k^2 M(p_j) + s_k D(p_j) + K(p_j)\right),
	\end{align}
	 $M(p_j), D(p_j), K(p_j) \in \mathcal{R}^{n \times n}$ and $C(p_j)^T, \ B(p_j) \in \mathcal{R}^n$. The variables $s_k$ and $p_j$ are the frequency variables  and parameters, respectively.}

	\myblue{PBTMR \cite{sonN2017} is based upon balanced truncation theory and is used for model reduction of parametric first and second order linear dynamical systems.
	Here, the second order system is transformed to the first order as
	\begin{align}\label{eq1:fst-order-ds}
	\begin{split}
	\mathbb{E}(p_j) \dot{z}(t) & = \mathbb{A}(p_j)z(t) + \mathbb{B}(p_j) u(t),\\
	y(t) & = \mathbb{C}(p_j)z(t),
	\end{split}
	\end{align}
	where 
	\begin{align}\label{eq:lyp}
	\begin{split}
	\mathbb{E}(p_j) = \begin{bmatrix}
	-K(p_j) & 0 \\
	0 & M(p_j)
	\end{bmatrix}, \ \ \mathbb{A}(p_j) = \begin{bmatrix}
	0 & -K(p_j) \\
	-K(p_j) & -D(p_j)
	\end{bmatrix}, \ \  \mathbb{B}(p_j) = \begin{bmatrix}
	0 \\
	B(p_j)
	\end{bmatrix}, \\ \mathbb{C}(p_j) = \begin{bmatrix}
	C_1(p_j) & C_2(p_j)
	\end{bmatrix}, \ 
	\dot{z}(t)
	 = \begin{bmatrix}
	\dot{x}(t) \\
	\ddot{x}(t)
	\end{bmatrix} \ 
	z(t)
	= \begin{bmatrix}
	x(t) \\
	\dot{x}(t)
	\end{bmatrix}.
		\end{split}
	\end{align}
Like in the case of RPMOR, here also one needs to build a projection matrix $V$. This requires solving Lyapunov equation of the form below for $Z(p_j)$ \cite{sonN2017}.
\begin{align}\label{eq:lyp_gen}
\begin{split}
\mathbb{A}(p_j)Z(p_j)\mathbb{E}^T(p_j) + \mathbb{E}(p_j)Z(p_j)\mathbb{A}^T(p_j) = & \ \mathbb{B}(p_j)\mathbb{B}^T(p_j) \quad or \\ 
vec\left(\mathbb{A}(p_j)Z(p_j)\mathbb{E}^T(p_j) + \mathbb{E}(p_j)Z(p_j)\mathbb{A}^T(p_j)\right) =  & \ vec\left(\mathbb{B}(p_j)\mathbb{B}^T(p_j)\right),\\
\end{split}
\end{align}	
where $vec$ denotes vectorization of a matrix  into a column vector. The above equation can be rewritten as
\begin{align}\label{eq:lyp_eq_gen}
A(p_j)z(p_j) = & \ b(p_j),
\end{align} 
where $A(p_j) = - \mathbb{A}(p_j) \otimes \mathbb{E} (p_j) - \mathbb{E}(p_j) \otimes \mathbb{A}(p_j), \ z(p_j) = vec \left(Z(p_j)\right)$, $b(p_j) = vec \left(\mathbb{B}(p_j)\mathbb{B}^T(p_j)\right)$, and $\otimes$ denotes the standard Kronecker product. Finally, the $V$ matrix is obtained as follow:
\begin{align*}
V  = & \ \left[z(p_1), z(p_1), \cdots, z(p_j)\right].
\end{align*} 
		}

\subsection{Solving Sequences of Linear Systems}
All these algorithms require solving sequences of the linear systems, 
 which is a key computational bottleneck when using them for reducing large dynamical systems. All the three algorithms lead to linear system matrices being \textit{dependent on parameters and have a similar form}. 

 The linear systems arising in RPMOR \cite{lihong2014} have the form as follows:
\begin{subequations}
\begin{align}
\label{eq:seq_LS_1}
\begin{split}
A{(1)} x{(0)}   & =  B \quad and    \\
A{(1)}  \left[x{(1)} \ \cdots \ x{(w+1)}\right] & =  \begin{bmatrix}
A_1 \\ \vdots \\ A_{w+1}
\end{bmatrix} x{(0)},  
\end{split}
\end{align}
\hspace{9cm}  \vdots  \\
 \begin{align}\label{eq:seq_LS_2}
 \begin{split}
  A{(\ell)} x{(0)}  & = B \quad and  \\
A{(\ell)}  \left[x{(1)} \ \cdots \ x{(w+1)}\right] & = \begin{bmatrix}
A_1 \\ \vdots \\ A_{w+1}
\end{bmatrix} x{(0)},   
\end{split}
\end{align}
\end{subequations}
where $B$ is given in (\ref{eq:Linear_param_system}); $A_1, A_2, \ldots, A_{w+1}$ are given in (\ref{eq:affine_system}); and $A{(1)}, \ A{(2)}, \ \ldots, \ A{(\ell)}$ for $\ell = 1, \ldots, \mathfrak{z}$ are given in (\ref{eq:param_mat1}), (\ref{eq:param_mat2}) $\&$ (\ref{eq:param_matl}). 
 
\myblue{
To compute the transfer function $H(s_k,p_j)$ in PMOR-L \cite{Ionita2014}, one needs to solve sequences of linear systems as
\begin{align}\label{eq:PMOR-L-eq}
\begin{split}
&\begin{cases}
A(s_1,p_1) x{(11)}&= B(p_1), \\
A(s_1,p_2) x{(12)}&= B(p_2), \\
&\vdots \\
A(s_1,p_j) x{(1j)} &= B(p_j), \\
\end{cases}\\
&\begin{cases}
A(s_2,p_1) x{(21)} &= B(p_1), \\ 
A(s_2,p_2) x{(22)} &= B(p_2), \\
&\vdots \\
A(s_2,p_j) x{(2j)} &= B(p_j), \\
\end{cases}\\
&\begin{cases}
A(s_k,p_1) x{(k1)} &= B(p_1), \\
A(s_k,p_2) x{(k2)} &= B(p_2), \\
&\vdots\\
A(s_k,p_j) x{(kj)} &= B(p_j), \\
\end{cases}
\end{split}
\end{align}
where $A(s_k,p_j) \in \mathcal{R}^{n \times n}$, $B(p_j) \in \mathcal{R}^n$ for $k = 1,\ldots, v$ and $j = 1,\ldots, w$ and are given in (\ref{eq:loewner_eq}) -- (\ref{eq:loewner_eq2}). 
This gives us the transfer function $ H(s_k,p_j) = C(p_j)^T [x_{11} \ x_{12} \cdots x_{kj} ]$.}

\myblue{
Solving the Lyapunov equations in PBTMR \cite{sonN2017} gives rise to the sequence of the linear systems as follows 
\begin{align}\label{eq:PMTMR-eq}
\begin{cases}
\begin{split}
A(p_1) z(1) &= b(p_1),\\
A(p_2) z(2) &= b(p_2),\\
&\vdots \\
A(p_j) z(j) & = b(p_j),
\end{split}
\end{cases}
\end{align}  
where $A(p_j) \in \mathcal{R}^{n^2 \times n^2}$ and $ b(p_j) \in \mathcal{R}^{n^2}$ for $j = 1, \ldots, w$ and are given in (\ref{eq1:fst-order-ds}) -- (\ref{eq:lyp_eq_gen}).
 }

\myblue{
If the dimensions of $A{(\ell)}$ for $\ell = 1,\ldots, \mathfrak{z}$, $A(s_k,p_j)$ and $A(p_j)$ for $k = 1,\ldots, v, \ \text{and} \ j = 1,\ldots, w$, are very large, one should use iterative methods (instead of direct method) to solve the above linear systems since they scale well. The time complexity of  direct methods is $\mathcal{O}(n^3)$ whereas for iterative methods it is $\mathcal{O}(n\cdot nnz)$, where $n$ represents the number of unknowns and $nnz$ is the number of non-zeros in  system matrix.
In RPMOR \cite{lihong2014}, the right hand side vectors of the second equations in (\ref{eq:seq_LS_1}) and (\ref{eq:seq_LS_2}) are available together. Hence, one can easily solve these linear systems simultaneously. For this, we can use a block version of the relevant iterative method \cite{OLEARY1980293, SIMONCINI1996457,Parks2006}. 
}

Preconditioning is a technique commonly used to accelerate the performance of iterative methods. In the algorithms above, the linear system matrices 
change with the change in parameters, however, this change is small. 

Since computing a new preconditioner for every new linear system is expensive, we  propose a cheap preconditioner update that avoids this. Here, we compute a preconditioner for the  initial linear system very accurately, and from the next linear systems, we use this initial preconditioner along with a cheap update. 
\myblue{People have proposed this for quantum Monte Carlo (QMC) \cite{ahuja2011improved}, model order reduction of non-parametric first order linear dynamical systems \cite{grim2015reusing,Wyatt2012} and model order reduction of non-parametric second order linear dynamical systems \cite{Navneet2016}, but not for parametric model order reduction (and hence, not specifically for model order reduction of parametric second order linear dynamical systems, which is our focus).	
	 The main challenge in this approach (cheap update) is to generate the \textit{best} sequence of preconditioners corresponding to the parametric coefficient matrices.  
}

The main contributions of this paper are as follows: Section \ref{sec:precon} discusses the use of iterative methods and preconditioners in this context. We propose our cheap preconditioner update techniques here as well. To support our theory, numerical results are provided in Section \ref{sec:numerical_results}. Finally, we give conclusions and discuss future work in Section \ref{sec:con_fut}.
\section{Our Approach} \label{sec:precon}
{Here, we first discuss preconditioned iterative methods in Section \ref{sec:prec_itr}, Next, for accelerating the iterative method, we discuss preconditioners in Section \ref{sec:precond}. We propose the theory of cheap preconditioner updates in Section \ref{sec:prec-update}. Finally, we discuss an application of preconditioner updates to the earlier discussed parametric model order reduction algorithms in Section \ref{sec:app_cheap_update} 
\subsection{Iterative Methods}\label{sec:prec_itr}
For solving linear systems of equations, either direct methods or iterative methods are used.  If a linear system is of large size, as discussed earlier, iterative methods are preferred over direct methods  because the latter is too expensive in terms of both storage and operation. 

Krylov subspace based methods are very popular class of iterative methods \cite{Greenbaum1997, Saad2003,van2003iterative}. There are many types of Krylov subspace methods. Some commonly used ones are Conjugate Gradient (CG), Generalized Conjugate Residual Orthogonal (GCRO), Generalized Minimal Residual (GMRES), Minimum Residual (MINRES), and BiConjugate Gradient (BiCG) etc. \cite{Greenbaum1997, Saad2003,van2003iterative}. The choice of method is problem dependent.

As discussed in Section \ref{sec:Intro}, if the linear systems have multiple right hand sides (available together), 
 then one can solve such linear systems by block iterative methods. 
This concept was introduced for the first time with Conjugate Gradient (CG) method \cite{OLEARY1980293}. A similar study with GMRES was proposed in \cite{SIMONCINI1996457}. 
Here, we give a brief overview of block iterative methods. Let the linear systems with multiple right-hand sides is given as
\begin{align*}
\mathcal{A}\mathcal{X}=\mathcal{B},
\end{align*}
where $\mathcal{A} \in \mathcal{R}^{n \times n}$, $\mathcal{B} \in \mathcal{R}^{n \times d_I}$, $d_I << n$. Given $R_0$ (i.e. $R_0 = \mathcal{B} - \mathcal{A}\mathcal{X}_0$) and $\mathcal{X}_0$ as the initial residual and the initial solution, respectively, these methods build the block Krylov subspace $\mathbb{K}^{\jmath} \left(\mathcal{A},R\right) = \it{span}\left \{R_0, \mathcal{A}R_0, \mathcal{A}^2R_0, \ldots, \mathcal{A}^{\jmath-1}R_0 \right \}$, and find solution in it \cite{OLEARY1980293,SIMONCINI1996457}.  

\subsection{Preconditioners}
\label{sec:precond}
Preconditioning is  used to accelerate the performance of iterative methods. If $P$ is a non-singular matrix that approximates the inverse of $\mathcal{A}$, that is $P \approx \mathcal{A}^{-1}$, then the system $\mathcal{A}P\tilde{\mathcal{X}} = \mathcal{B} \ \text{with}  \ \mathcal{X} = P \tilde{\mathcal{X}}$ may be faster to solve than the original one (i.e. $\mathcal{A}\mathcal{X} = \mathcal{B}$)\footnote[1]{Here, we use right preconditioning, i.e. the preconditioner is applied to the right of the linear system matrix. Similar analysis can be done with left preconditioning, i.e. with the preconditioner on the left of the linear system matrix.}. For most of the dynamical systems reduced by the earlier discussed algorithms \cite{lihong2014,Ionita2014,sonN2017}, the iterative methods stagnate or are slow in convergence  (see Numerical Results section). Hence, we use a preconditioner.

Besides making the system easier to solve by an iterative method, a preconditioner should be cheap to construct and apply. Some existing preconditioning techniques include Successive Over Relaxation, Polynomial Based, Incomplete Factorizations, Sparse Approximate Inverse (SPAI), and Algebraic Multi-Grid \cite{Benzi2002}.  
 SPAI preconditioners are known to work in the most general setting and can be easily parallelized \cite{Benzi2002,saadchow}. Among the others, Incomplete Factorizations are also general but these cannot be easily parallelized. Hence, we use a parallel version of SPAI.  We briefly discuss SPAI preconditioner next. 

In constructing a preconditioner $P$ for a coefficient matrix $\mathcal{A}$, we would like $ \mathcal{A}P \approx I$ ($I$ is the identity matrix). SPAI preconditioner finds $P$ by minimizing the associated error norm $\|I-\mathcal{A}P\|$. If the norm used is Frobenius norm, then the minimization problem becomes
\begin{align*}
\min_{P} \|I-\mathcal{A} P\|_f.
\end{align*}   
This minimization problem can be rewritten as  
\begin{align*}
\min_{P} \|I-\mathcal{A} P \|_f^2 = \min_{p^\imath} \sum_{\imath = 1}^{n} \|e^\imath-\mathcal{A} p^\imath\|_2^2,
\end{align*}
where $p^\imath$ and $e^\imath$ are the $\imath^{th}$ columns of $P$ and $I$ matrices, respectively. This minimization problem is just a least square problem, to be solved for $n$ different right hand sides \cite{saadchow,Alexander2007}.

\subsection{Theory of Cheap Preconditioner Updates}
\label{sec:prec-update}
\myblue{ In general the sequences of linear systems (\ref{eq:seq_LS_1}-\ref{eq:seq_LS_2}), (\ref{eq:PMOR-L-eq}) and (\ref{eq:PMTMR-eq}) can be written as
	\begin{align}\label{eq:general_sq}
	\begin{split}
\mathcal{A}{(1)} \mathcal{X}{(1)} & = \mathcal{B}{(1)},\\
\mathcal{A}{(2)} \mathcal{X}{(2)}&  = \mathcal{B}{(2)},\\
& \vdots\\
\mathcal{A}{(i)} \mathcal{X}{(i)} & = \mathcal{B}{(i)},
	\end{split}
	\end{align}
where $\mathcal{A}{(1)}, \ldots, \mathcal{A}{(i)} \in \mathcal{R}^{n \times n}$ and $\mathcal{B}{(1)}, \ldots, \mathcal{B}{(i)} \in \mathcal{R}^n$. Let $P_1$ be a good preconditioner for $\mathcal{A}(1)$ (i.e. computed by $\underset{P_1}{min}\|I- \mathcal{A}(1)P_1\|$) then, $P_i$ (preconditioner corresponding to $\mathcal{A}(i)$, for $i = 2, \ldots, m$) can be computed as given in Table \ref{tab:cheap_update}.} 
\begin{table}[]
	\centering
	\caption{Cheap Preconditioner Update Approaches }
	\def\arraystretch{1.3}
	\begin{tabular}{|l|l|}
		\hline
		First Approach & Second Approach \\ \hline
		\textbullet \ \ $\mathcal{A}(1)P_1 = \mathcal{A}(2)P_2$ &  \\
		\textbullet \ \ If $P_2 = Q_2P_1$, & \qquad Same as the \textit{first} approach\\
		\qquad then $ \mathcal{A}(1)P_1 = \mathcal{A}(2)Q_2P_1$ & \\
		\textbullet \ \ $\min\limits_{Q_2} \|\mathcal{A}(1)-\mathcal{A}(2) Q_2\|_f^2$ &   \\
		 \hline
		 \textbullet \ \ $\mathcal{A}(1)P_1 = \mathcal{A}(3)P_3$ & \textbullet \ \ $\mathcal{A}(2)P_2 = \mathcal{A}(3)P_3$ \\
		 \textbullet \ \ If $P_3 = Q_3 P_1  $, & \textbullet \ \ If $P_3 = Q_3 P_{2}$, \\
		 \qquad then $ \mathcal{A}(1)P_1 = \mathcal{A}(3)Q_3P_1$ & \qquad then $\mathcal{A}(2)P_2 = \mathcal{A}(3)Q_3P_2$ \\
		\textbullet \ \ $\min\limits_{Q_3} \| \mathcal{A}(1)- \mathcal{A}(3) Q_3\|_f^2$ & \textbullet \ \ $\min\limits_{Q_3} \|\mathcal{A}{(2)}-\mathcal{A}(3)Q_3\|_f^2$  \\
		 \hline		
		\qquad \qquad \qquad \vdots & \qquad \qquad \qquad \vdots \\ \hline
		\textbullet \ \ $\mathcal{A}(1)P_1 = \mathcal{A}(i)P_{i}$ & \textbullet \ \ $\mathcal{A}(i-1)P_{i-1} = \mathcal{A}(i)P_{i}$ \\
		\textbullet \ \ If $P_{i} = Q_{i} P_1  $, & \textbullet \ \ If $P_{i} = Q_{i} P_{i-1} $, \\
		\qquad then $ \mathcal{A}(1)P_1 = \mathcal{A}(i)Q_{i}P_1$ & \qquad then $\mathcal{A}(i-1)P_{i-1} = \mathcal{A}(i)Q_{i}P_{i-1}$ \\
		\textbullet \ \ $\min\limits_{Q_{i}} \|\mathcal{A}(1) - \mathcal{A}(i) Q_{i}\|_f^2$& \textbullet \ \ $\min\limits_{Q_{i}} \|\mathcal{A}{(i-1)}-\mathcal{A}(i)Q_{i}\|_f^2$  \\
		\hline	
	\end{tabular}
	\label{tab:cheap_update}
\end{table}

\myblue{The approaches provided in Table \ref{tab:cheap_update} have competing trade-offs. In the \textit{first} approach, the minimization is harder to solve because  $\mathcal{A}{(i)}$ and $\mathcal{A}{(i-1)}$ would be closer than $\mathcal{A}(i)$ and $\mathcal{A}(1)$ (since the sequences of matrices in (\ref{eq:general_sq}) change slowly). However, $P_{i}$ is more accurate since $P_1$ is very accurate  (see $\min\limits_{P_1}  \|I- \mathcal{A}(1)P_1\|$ ). 
In the \textit{second} approach, the minimization is easier using the same argument as earlier, while the preconditioner at each step is less accurate ($P_{i}$ is formed from $P_{i-1}$, which already has approximation errors). 
}

\myblue{ In non-parametric model order reduction, the relative difference between $\mathcal{A}(i)$ and $\mathcal{A}{(i-1)}$ is substantial, which means $\mathcal{A}(1)$ and $\mathcal{A}(i)$ are further away (expansion points change rapidly, however, they are still ``close'' to be able to apply cheap update (see \cite{Wyatt2012,Navneet2016})). Hence, the \textit{first} approach is not very efficient there and the \textit{second} approach fits well \cite{Navneet2016}. In the parametric case, change in $\mathcal{A}{(i-1)}$ to $\mathcal{A}(i)$ is such that $\mathcal{A}(1)$ and $\mathcal{A}(i)$ are not too far (or the relative difference between $\mathcal{A}(1)$ and $\mathcal{A}(i)$ is not substantial). Hence, the minimization problem $\min\limits_{Q_{i}} \|\mathcal{A}(1)- \mathcal{A}(i)Q_{i}\|$ is almost as easy to solve  as $\min\limits_{Q_{i}} \| \mathcal{A}{(i-1)}- \mathcal{A}(i)Q_{i}\|$. Since, the \textit{first} approach has an extra advantage of less loss of accuracy during building the preconditioner after minimization ($P_{i} = Q_{i}P_1$ as compared to $ P_{i} = Q_{i} P_{i-1}$), we propose its use. The experimental results support our this argument as well. 
}

To summarize, when using basic SPAI we need to solve $\min\limits_{P_{i}} \|I- \mathcal{A}{(i)}P_i\|^2$, which we transform to first $\min\limits_{P_i} \|\mathcal{A}{(1)}P_1 - \mathcal{A}{(i)}P_i\|^2$, and subsequently to $\min\limits_{{Q}_i}\|\mathcal{A}{(1)}-\mathcal{A}{(i)} {Q}_i\|^2$. This last formulation is usually much easier to solve since $\mathcal{A}{(1)}$ and $\mathcal{A}{(i)}$ are close to each other (change in expansion points only), as compared to the first formulation where $I$ and $\mathcal{A}{(i)}$ could be very different. 

\subsection{Application of Cheap Preconditioner Updates} \label{sec:app_cheap_update}

Recall (\ref{eq:seq_LS_1})-(\ref{eq:seq_LS_2}) in RPMOR \cite{lihong2014}, let $A{(1)} = A_0 + \tilde{s}_1^{1} A_1 + \tilde{s}_2^1 A_2 + \ldots + \tilde{s}_{w+1}^1 A_{w+1}$ and $A{(\ell)} = A_0 + \tilde{s}_1^{\ell} A_1 + \tilde{s}_2^{\ell} A_2 + \ldots + \tilde{s}_{w+1}^{\ell} A_{w+1}$ be  two coefficient matrices for different expansion points $\tilde{\tilde{s}}_1 =\left[\tilde{s}_1^1, \cdots, \tilde{s}_{w+1}^1\right]$ and $\tilde{\tilde{s}}_{\ell} = \left[\tilde{s}_1^{\ell}, \cdots, \tilde{s}_{w+1}^{\ell}\right]$, respectively. If the difference between $A{(1)}$  and $A{(\ell)}$ is small, then one can exploit this while building preconditioners for this sequence of matrices. 

Let $P_1$ be a good initial  preconditioner for $A{(1)}$. 
Then, a cheap preconditioner update can be obtained  by making $A{(1)}P_1 \approx A{(\ell)}P_{\ell}$, where $\ell \in \{1,\ldots, \mathfrak{z}\}$, and $\mathfrak{z}$ denotes the number of expansion points.    
Expressing $A{(\ell)}$ in terms of $A{(1)}$ we get
\begin{equation*}
	\begin{split}
	A{(\ell)}  =  & A{(1)}  \left[ I+ \left(\tilde{s}_1^{\ell} - \tilde{s}_1^1\right)\left(A{(1)}\right)^{-1}A_1+\left(\tilde{s}_2^{\ell} - \tilde{s}_2^1\right)\left(A{(1)}\right)^{-1}  A_2  + \cdots + \right. 
    \left.   \left(\tilde{s}_{w+1}^{\ell} - \tilde{s}_{w+1}^1\right)\left(A{(1)}\right)^{-1}A_{w+1} \right].
		\end{split}
\end{equation*}
Now we enforce $A{(1)} P_1 = A{(\ell)} P_{\ell}$ or
\begin{equation*}
\begin{split}
A{(1)} P_1   =  & \ A{(1)}  \left[ I+ \left(\tilde{s}_1^{\ell} - \tilde{s}_1^1\right)\left(A{(1)}\right)^{-1}A_1 
+ \left(\tilde{s}_2^{\ell} - \tilde{s}_2^1\right)\left(A{(1)}\right)^{-1}  A_2  + \cdots +      \left(\tilde{s}_{w+1}^{\ell} - \tilde{s}_{w+1}^1\right) \right.  \left. \left(A{(1)}\right)^{-1}A_{p+1} \right] \cdot \\
 &\left[ I+ \left(\tilde{s}_1^{\ell} - \tilde{s}_1^1\right)  \left(A{(1)}\right)^{-1}A_1   
   +\left(\tilde{s}_2^{\ell} - \tilde{s}_2^1\right)\left(A{(1)}\right)^{-1}A_2   + \cdots + \right. \left.
    \left(\tilde{s}_{w+1}^{\ell} - \tilde{s}_{w+1}^1\right)\left(A{(1)}\right)^{-1}A_{w+1}\right]^{-1}  P_1, \\
  = & \ A{(\ell)} P_{\ell},
\end{split}
\end{equation*}
where  
\begin{equation*}
\begin{split}
  P_{\ell} = & \left[ I+ \left(\tilde{s}_1^{\ell} - \tilde{s}_1^1\right)  \left(A{(1)}\right)^{-1}A_1   
  +\left(\tilde{s}_2^{\ell} - \tilde{s}_2^1\right)\left(A{(1)}\right)^{-1}A_2   + \cdots + \right. \left.
  \left(\tilde{s}_{w+1}^{\ell} - \tilde{s}_{w+1}^1\right)\left(A{(1)}\right)^{-1}A_{w+1} \right]^{-1}  P_1. 
 \end{split}
 \end{equation*}
 
Let 
\begin{equation*}
\begin{split}
{Q}_{\ell} = & \left[ I+ \left(\tilde{s}_1^{\ell} - \tilde{s}_1^1\right)  \left(A{(1)}\right)^{-1}A_1   
+\left(\tilde{s}_2^{\ell} - \tilde{s}_2^1\right)\left(A{(1)}\right)^{-1}A_2   + \cdots + \right. \left.
\left(\tilde{s}_{w+1}^{\ell} - \tilde{s}_{w+1}^1\right)\left(A{(1)}\right)^{-1}A_{w+1} \right]^{-1},
  \end{split}
  \end{equation*} 
   then the above implies $A{(1)} P_1= A{(\ell)} {Q}_{\ell}P_1$.
   This leads us to the idea that instead of solving  $A{(1)} P_1 \approx A{(\ell)} {Q}_{\ell} P_1$ for ${Q}_{\ell}$ leading to $P_{\ell} = {Q}_{\ell}P_1$, we solve a simpler problem given below
\begin{align*}
\min_{{Q}_{\ell}}\left|\left|A{(1)} - A{(\ell)} {Q}_{\ell}\right|\right|^{2}_{f}=\min_{\left(q_{\ell}\right)^{(\imath)}} \sum_{\imath=1}^{n}\left|\left|\left(a_1\right)^{(\imath)}-A{(\ell)} \left(q_{\ell}\right)^{(\imath)}\right|\right|^{2}_2,
\end{align*}
 where $\left(a_1\right)^{(\imath)}$ and $\left(q_{\ell}\right)^{(\imath)}$ denote the $\imath^{th}$ columns of $A{(1)}$ and ${Q}_\ell$, respectively.

\myblue{
Next, we look at a cheap update for PMOR-L \cite{Ionita2014}. Here, we can express the relation between coefficient matrices of the two consecutive linear systems by the following two ways: First, by capturing the changes in frequency variables $(s_k)$, and second, by capturing the changes in parameters $(p_j)$. In general, savings in computation time is less in the first case. This is because, as earlier, frequency variables change more rapidly than parameters. We discuss both cases here in the above order.  
}      

\myblue{
Let $ A(s_1,p_1) $ and $A(s_k,p_1)$ be two coefficient matrices for parameter $p_1$ and different values of frequency parameters  $s_k$ for $k = 1,\ldots, v$ (recall (\ref{eq:loewner_eq2}) and (\ref{eq:PMOR-L-eq})). Now, expressing $ A(s_k,p_1) $ in terms of $ A(s_1,p_1) $ 
we get
\begin{equation*}
\begin{split}
A(s_k,p_1)  =  & \ A(s_1,p_1)  \left[ I+ (s_k^2-s_1^2)A(s_1,p_1)^{-1}M(p_1)  + (s_k-s_1)A(s_1,p_1)^{-1}D(p_1)\right].
\end{split}
\end{equation*}
Enforcing $A(s_1,p_1) P_{11} = A(s_k,p_1) P_{k1}$ we have
\begin{align*}
A(s_1,p_1) P_{11} = &  \ A(s_1,p_1) \left[ I+ (s_k^2-s_1^2)A(s_1,p_1)^{-1}M(p_1)  + (s_k-s_1)A(s_1,p_1)^{-1}D(p_1)\right] \\
&  \left[ I+ (s_k^2-s_1^2)A(s_1,p_1)^{-1}M(p_1)  + (s_k-s_1)A(s_1,p_1)^{-1}D(p_1)\right]^{-1} P_{11}, \\
= & \ A(s_k,p_1) P_{k1},
\end{align*}
where $P_{k1} = \left[ I+ (s_k^2-s_1^2)A(s_1,p_1)^{-1}M(p_1)  + (s_k-s_1)A(s_1,p_1)^{-1}D(p_1)\right]^{-1}P_{11}$.\\ 
Let $Q_{k1} =   \left[ I+ (s_k^2-s_1^2)A(s_1,p_1)^{-1}M(p_1)  + (s_k-s_1)A(s_1,p_1)^{-1}D(p_1)\right]^{-1}$, then the above implies $A{(s_1,p_1)} P_{11}= A{(s_k,p_1)} {Q}_{k1}P_{11}$.
This leads us to the idea that instead of solving  $A{(s_1,p_1)} P_{11} \approx A{(s_k,p_1)} {Q}_{k1}P_{11}$ for ${Q}_{k1}$ leading to $P_{k1} = {Q}_{k1}P_{11}$, we solve a simpler problem given below
\begin{align*}
\min_{Q_{k1}}\|A(s_1,p_1) - A(s_k,p_1)Q_{k1}\|_f^2.
\end{align*}
Similarly, for any parameter $p_j$ $(\text{for}\ j = 2, \ldots, w)$ we solve for $Q_{kj}$ from 
  \begin{align*}
  \min_{Q_{kj}}\|A(s_1,p_j) - A(s_k,p_j)Q_{kj}\|_f^2,
  \end{align*}
  leading to $P_{kj} = Q_{kj} P_{1j}$.
}

\myblue{In the second case, we are unable to express the two linear systems in-terms of each other, i.e. $A(s_1,p_1)$ in-terms of $A(s_1,p_j)$ or $A(s_k,p_1)$ in-terms of $A(s_k,p_j)$ unless we know how the matrices $M$ and $D$ from (\ref{eq:Linear_param_system}) depend on $p_j$. However, this can be easily worked out once the input dynamical system is known. We give this derivation for a commonly used example from \cite{Ionita2014} (the paper which proposed PMOR-L) in Appendix I.
}

\myblue{ A cheap update for PBTMR \cite{sonN2017} can be worked out as below. Let $ A(p_1) $ and $A(p_j) $ be two coefficient matrices for $j = 2, \ldots, w$ $\left( \text{recall} \ (\ref{eq:lyp}), (\ref{eq:lyp_gen}) \ \text{and} \ (\ref{eq:PMTMR-eq}) \right)$ 
\begin{align*}
 A(p_1) & = -\mathbb{E}(p_1)  \otimes \mathbb{A}(p_1) - \mathbb{A}(p_1) \otimes \mathbb{E}(p_1), \\
					& =	- \begin{bmatrix}
					-K(p_1) & 0 \\
					0 & M(p_1)
					\end{bmatrix} \otimes \begin{bmatrix}
					0 & -K(p_1) \\
					-K(p_1) & -D(p_1)
					\end{bmatrix} - \begin{bmatrix}
					0 & -K(p_1) \\
					-K(p_1) & -D(p_1)
					\end{bmatrix} \otimes  	\begin{bmatrix}
					-K(p_1) & 0 \\
					0 & M(p_1)
					\end{bmatrix},\\
					& = \begin{bmatrix}
					0 & -K^2(p_1) & -K^2(p_1) & 0 \\
					-K^2(p_1) & -K(p_1)D(p_1) & 0 & K(p_1)M(p_1) \\
					-K^2(p_1) & 0 & -D(p_1)K(p_1) & M(p_1)K(p_1) \\
					0 & K(p_1)M(p_1) & M(p_1)K(p_1) &  M(p_1)D(p_1) + D(p_1)M(p_1)
					\end{bmatrix},\\
A(p_j) & = -\mathbb{E}(p_j)  \otimes \mathbb{A}(p_j) - \mathbb{A}(p_i) \otimes \mathbb{E}(p_j), \\
& =	- \begin{bmatrix}
-K(p_j) & 0 \\
0 & M(p_j)
\end{bmatrix} \otimes \begin{bmatrix}
0 & -K(p_j) \\
-K(p_j) & -D(p_j)
\end{bmatrix} - \begin{bmatrix}
0 & -K(p_j) \\
-K(p_j) & -D(p_j)
\end{bmatrix} \otimes  	\begin{bmatrix}
-K(p_j) & 0 \\
0 & M(p_j)
\end{bmatrix},\\
& = \begin{bmatrix}
0 & -K^2(p_j) & -K^2(p_j) & 0 \\
-K^2(p_j) & -K(p_j)D(p_j) & 0 & K(p_j)M(p_j) \\
-K^2(p_j) & 0 & -D(p_j)K(p_j) & M(p_j)K(p_j) \\
0 & K(p_j)M(p_j) & M(p_j)K(p_j) &  M(p_j)D(p_j)+D(p_j)M(p_j)
\end{bmatrix}.					
\end{align*}
Here, we are unable to express the $j^{th}$ linear system in-terms of the first linear system (i.e. $A(p_j)$ in-terms of $A(1)$). However, we can see that the two  linear systems have structural similarities (only parameters are varying). 
We can abstract out the relationship between these linear systems if structure of $\mathbb{E}(p_j)$ and $\mathbb{A}(p_j)$ is more 
explicitly known. 
 Hence, we have worked out a cheap preconditioner update for a popular example as used in the paper that proposed PBTMR, i.e. \cite{sonN2017}, in Appendix II. 
}
 

\section{Numerical Results}
\label{sec:numerical_results}

We demonstrate our proposed preconditioned iterative solver theory using RPMOR \cite{lihong2014} as our candidate parametric model order reduction algorithm and a micro-gyroscope model \cite{lihong2013} as our test dynamical system (earlier version of RPMOR in \cite{lihong2013} is tested on a micro-gyroscope model). This model is a parametric Single Input Single Output (SISO) second order linear dynamical systems of size $17, 931$, and is given as
\begin{equation*}
\begin{split}
s^2 M(d)x + s D(\theta, \alpha, \beta, d) x + K(d) x = B u(s), \\
y = Cx, 
\end{split}
\end{equation*}
where $M(d) = M_1 + d M_2, \ D(\theta, \alpha, \beta, d) = \theta (D_1 + d D_2) + \alpha M(d) + \beta K(d), \ K(d) = K_1 + (1/d)K_2 + d K_3.$ In above equation, there are eleven  variables and all must be considered as individual parameters.
These  parameters at the ${\ell}^{th}$ expansion point are $\tilde{s}_1^{\ell} = s^2, \tilde{s}_2^{\ell} = s^2 d, \tilde{s}_3^{\ell} = s\theta, \tilde{s}_4^{\ell} = s \theta d, \tilde{s}_5^{\ell} = s \alpha, \tilde{s}_6^{\ell} = s \alpha d, \tilde{s}_7^{\ell} = s \beta, \tilde{s}_8^{\ell} = s\beta/d, \tilde{s}_9^{\ell} = s\beta d, \tilde{s}_{10}^{\ell} = 1/d $ and $\tilde{s}_{11}^{\ell} = d$. Usually $\alpha$ and $\beta$ are taken as zero \cite{lihong2013}, and hence we are left with six parameters $\tilde{s}_1^{\ell}, \tilde{s}_2^{\ell}, \tilde{s}_3^{\ell}, \tilde{s}_4^{\ell}, \tilde{s}_{10}^{\ell}$ and $\tilde{s}_{11}^{\ell}$, where $\ell = 1, 2, \ldots, \mathfrak{z}$.

We reduce this model to size $304$, and use four expansion points (i.e. $\mathfrak{z} = 4$) $\tilde{\tilde{s}}_{\ell} = \left[\tilde{s}_1^{\ell}, \tilde{s}_2^{\ell}, \tilde{s}_3^{\ell}, \tilde{s}_4^{\ell}, \tilde{s}_{10}^{\ell}, \tilde{s}_{11}^{\ell}\right]$ for $\ell = 1, 2, \ldots, 4$ based upon values in \cite{lihong2013}.
That is,
\\ $\tilde{\tilde{s}}_1 = \left[-4\pi^2 \times 0.065, -4\pi^2 \times 0.065, 5 \pi \sqrt{-1} \times 10^{-7}, 5 \pi \sqrt{-1}  \times 10^{-7},1, 1 \right],\\ \tilde{\tilde{s}}_2 = \left[- 4\pi^2 \times 0.065, - 8 \pi^2 \times 0.065, 5 \pi \sqrt{-1} \times 10^{-7}, 10 \pi \sqrt{-1} \times 10^{-7}, 0.5, 2 \right],\\ \tilde{\tilde{s}}_3 = \left[- 4 \pi^2 \times 0.0225, - 8 \pi^2 \times 0.0225, 3 \pi \sqrt{-1}  \times 10^{-7}, 6 \pi \sqrt{-1} \times 10^{-7}, 0.5, 2 \right],$ and \\$\tilde{\tilde{s}}_4 = \left[- 4\pi^2 \times 0.0225, - 4 \pi^2 \times 0.0337, 3 \pi \sqrt{-1}  \times 10^{-7}, 4.5 \pi \sqrt{-1} \times 10^{-7}, 0.66, 1.5 \right]$. 

The linear systems that arise here are of size $17,931 \times 17,931$ with non-symmetric linear system matrices. As earlier, we use iterative methods instead of direct methods. 
Of the many available iterative methods for solving non-symmetric linear systems, we use GCRO \cite{DESTURLER1996}. 
In fact, we use block GCRO \cite{parks2011,Michael2005} because of the reasons discussed before (availability of the multiple right hand sides together). 
We also compare usage of GCRO and block GCRO. The stopping tolerance is taken as $10^{-10}$ for all cases.

{Preconditioning has to be employed when iterative methods fail or have a slow convergence. Here, for this model, we observe that unpreconditioned GCRO fails to converge. 
We use a Modified Sparse Approximate Inverse (MSPAI 1.0) proposed in \cite{Alexander2007} as our SPAI preconditioner. This is because MSPAI uses a linear algebra library for solving sparse least squares problems, which arise here. 
We use standard initial settings of MSPAI as follows: tolerance (ep) of $0.0001$ and cache size (cs) of 80.}  



{We test on a machine with the following configuration: Intel Xeon (R) CPU E5-1620 V3 $@$ 3.50 GHz., frequency 1200 MHz., 8 CPU and 64 GB RAM.   
All the codes are written in MATLAB (2016b) (including RPMOR, GCRO and block GCRO) except SPAI and SPAI Update. MATLAB is used because of ease of rapid prototyping. Computing SPAI and SPAI update in MATLAB is expensive, therefore, we use C++ version of these (SPAI from MSPAI and SPAI update written by us). MSPAI further uses BLAS, LAPACK and ATLAS libraries. While solving every linear system in the sequence, we first compute SPAI and SPAI update separately and save them.  Then, we run MATLAB code along with the saved preconditioner matrices (i.e. SPAI and SPAI update).} 



\subsection{Analysis} \label{sec:anal}
Here, we quantify the difference between the parameterized coefficient matrices because of changing expansion points.   
The first coefficient matrix,  $A{(1)}$ defined in (\ref{eq:param_mat1}), is given as   
\begin{align*}
A{(1)} = K_1 + s_1^1 M_1 + s_2^1 M_2 + s_3^1 D_1 + s_4^1 D_2  + s_{10}^1 K_2 + s_{11}^1 K_3,
\end{align*}
and other coefficient matrices, $A{(\ell)}$ for $\ell = 2, \ldots, 4$, are given as
\begin{align*}
A{(\ell)} = K_1 + s_1^{\ell} M_1 + s_2^{\ell} M_2 + s_3^{\ell} D_1 + s_4^{\ell} D_2  + s_{10}^{\ell} K_2 + s_{11}^{\ell} K_3.
\end{align*}

{Now, we analyze how SPAI update is beneficial. As discussed in Section \ref{sec:prec-update}, SPAI update is useful when $\|I - A{(\ell)}\|_f$ is large and $\|A{(1)} - A{(\ell)}\|_f$ is small. This data for the above four expansion points is given in Table \ref{tab:spai_update_analysis}. From column 3, it is observed that the change in $\|I - A{(\ell)}\|_f$, for $\ell = 1, 2,\ldots, 4$, is large, whereas from column 4, we can see that the change  from one expansion point  to another for $\|A{(1)} - A{(\ell)}\|_f$ is small. 
	}    
\subsection{Iteration Count and Computation Time Comparison} \label{sec:itr_comp}

{First, we compare GCRO and block GCRO method in Table \ref{Tab:GCRO_Block_GCRO} (in-terms of iteration count and computation time).  
Here, we have to solve $43$ linear systems at each step in the RPMOR algorithm. That is, GCRO is executed $43$ times at each RPMOR algorithm step, while block GCRO is executed  only 8 times (single linear system is solved in the first call to block GCRO; further $6$ linear systems are solved together for each of remaining calls to block GCRO). 
Iteration counts for both the solvers (GCRO and block GCRO) are given in columns 3 and 6, respectively. Computation times corresponding to these solvers are given in columns 4 and 7, respectively. In this table, we give iteration count as well as computation time for solving each linear system (on an average). This is because we can extract the needful information from average also. 
From the last row of Table \ref{Tab:GCRO_Block_GCRO} it is clear that block GCRO saves nearly $80 \%$ in both iteration count as well as computation time as  compared to GCRO.}

{Next, in Table \ref{tab:SPAI_CT}, we provide the computation time of  SPAI and SPAI update\footnote{The underlying linear solver, GCRO or block GCRO, does not affect this.}. In the first column, we give the step number of the RPMOR algorithm that corresponds to the four expansion points.  The computation time of SPAI  and SPAI update is given in columns 2 and 3, respectively, corresponding to each RPMOR step. 
Here, both SPAI and SPAI update are computed once at each RPMOR step  and applied to all linear systems (43 linear systems) because the coefficient matrices are not changing (only right hand sides change).} 
{At the first RPMOR step, SPAI and SPAI update take the same amount of computation time. This is because SPAI update is not applicable at this step. From second RPMOR step onwards, we see substantial savings with SPAI update as compared to SPAI (approximately $70 \%$).}

{Table \ref{tab:bGCRO_DR_SPAI} gives the computation time of block GCRO with SPAI and block GCRO with SPAI update. In the first column we give the step number of the RPMOR algorithm  corresponding to the four expansion points. The computation times of block GCRO with SPAI and block GCRO with SPAI update are given in columns 2 and 3, respectively. As above, at the first RPMOR step, the computation time of block GCRO with SPAI and block GCRO with SPAI update is the same. From the second step onwards, we see substantial savings in computation time by using block GCRO with SPAI update over block GCRO with SPAI. From the last row of this table, it is clear that block GCRO with SPAI update saves approximately $65 \%$ computation time as compared to block GCRO with SPAI.}
\begin{table}[]
	\centering
	\caption{SPAI and SPAI Update Analysis}
	\label{tab:spai_update_analysis}
	\begin{tabular}{|c|c|c|c|}
		\hline
		\begin{tabular}[c]{@{}c@{}}RPMOR  \\ Steps\end{tabular} & Expansion Points & $\|I-A{(\ell)}\|_f$ & $\|A{(1)}-A{(\ell)}\|_f$ \\ \hline
		1 & $\tilde{\tilde{s}}_1$ & $1.36 \times 10^{02}$ & $0$ \\ \hline
		2 & $\tilde{\tilde{s}}_2$ & $1.36 \times 10^{02}$ & $7.9 \times 10^{-04}$ \\ \hline
		3 & $\tilde{\tilde{s}}_3$ & $1.36 \times 10^{02}$ & $7.9 \times10^{-04}$ \\ \hline
		4 & $\tilde{\tilde{s}}_4$ & $1.36 \times 10^{02}$ & $4.9 \times 10^{-04}$ \\ \hline
	\end{tabular}
\end{table}
\begin{table}[]
	\centering
	\caption{GCRO and block GCRO Iteration Count and Computation Time}
	\label{Tab:GCRO_Block_GCRO}
	\begin{tabular}{|c|c|c|c|c|c|c|}
		\hline
		\multirow{2}{*}{\begin{tabular}[c]{@{}c@{}}RPMOR\\ Steps\end{tabular}}  & \multicolumn{3}{c|} {GCRO}   & \multicolumn{3}{c|}{block GCRO} \\ \cline{2-7}  
		
		& \begin{tabular}[c]{@{}c@{}}No. of iterative  \\ solver calls \end{tabular} & \multicolumn{1}{c|}{\begin{tabular}[c]{@{}c@{}} Iteration \\ Count$^\dagger$ \end{tabular}} & \begin{tabular}[c]{@{}c@{}} Computation \\ Time$^\dagger$ \\ (seconds)\end{tabular} & \begin{tabular}[c]{@{}c@{}} No. of iterative \\ solver calls$^\ddagger$ \\ \end{tabular} & \begin{tabular}[c]{@{}c@{}} Iteration \\ Count$^\dagger$   \end{tabular} & \begin{tabular}[c]{@{}c@{}} Computation \\ Time$^\dagger$ \\(seconds)\end{tabular} \\ \hline
		1 & 43 & 426 & 14 & 08  & 528 & 15 \\ \hline
		2 & 43 & 404 & 14 & 08  & 502 & 15 \\ \hline
		3 & 43 & 406 & 14 & 08 & 504 &  15 \\ \hline
		4 & 43 & 412 & 15 & 08 & 508  & 15 \\ \hline
		Sum & 172 & 1648 & 57 & 32 & 2042  &  60  \\  
		\specialrule{0.1em}{0.1em}{0.1em}
      Total & $172$ & \begin{tabular}[c]{@{}c@{}} $43 \times 1648$ \\ $= 70864$ \end{tabular} & \begin{tabular}[c]{@{}c@{}} $43 \times 57$  \\ $= 2451$ or \\ 41 (minutes)  \end{tabular} & $32$ &  \begin{tabular}[c]{@{}c@{}} $08 \times 2042$  \\ $= 16336$ \end{tabular}  &  \begin{tabular}[c]{@{}c@{}} $08 \times 60$  \\ $= 480$ or \\ 08 (minutes)  \end{tabular}\\ \hline
	\end{tabular}
	\begin{tablenotes}
		\item[1] $\dagger$ per iterative solve on an average.
		\item[2] $\ddagger$ single linear system is solved in the first call to block GCRO; further 6 linear systems are solved together for each of the remaining calls to block GCRO.
	\end{tablenotes}
\end{table}

\begin{table}[!]
	\centering
	\setlength\tabcolsep{4pt}
	\begin{minipage}{0.5\textwidth}
		\centering
		\caption{Computation Time of SPAI and SPAI update}
		\label{tab:SPAI_CT}
		\begin{tabular}{|c|c|c|}
			\hline
			\multirow{2}{*}{\begin {tabular}[c]{@{}c@{}}RPMOR\\ Steps\end{tabular}} & \multicolumn{2}{c|}{\begin{tabular}[c]{@{}c@{}}Computation Time \\ (minutes)\end{tabular}} \\ \cline{2-3} 
		& SPAI & SPAI Update \\ \hline
		1 & 27 & 27 \\ \hline
		2 & 26 & 2 \\ \hline
		3 & 26 & 2 \\ \hline
		4 & 26 & 2 \\ \hline
		Total Time & 105 & 33 \\ \hline
	\end{tabular}
	\end{minipage}%
	\hfill
	\begin{minipage}{0.5\textwidth}
		\centering
		\caption{Computation Time of block GCRO with SPAI and  SPAI Update} 
		\label{tab:bGCRO_DR_SPAI} 
		\begin{tabular}{|c|c|c|}
			\hline
			\multirow{2}{*}{\begin{tabular}[c]{@{}c@{}}RPMOR\\ Steps\end{tabular}} & \multicolumn{2}{c|}{\begin{tabular}[c]{@{}c@{}}Computation Time \\ (minutes)\end{tabular}} \\ \cline{2-3} 
			& \begin{tabular}[c]{@{}c@{}} block GCRO \\ plus SPAI \end{tabular} & \begin{tabular}[c]{@{}c@{}}block GCRO \\ plus SPAI Update \end{tabular} \\ \hline
			1 & 29 & 29 \\ \hline
			2 & 28 & 4 \\ \hline
			3 & 28 & 4 \\ \hline
			4 & 28 & 4 \\ \hline
			Total Time & 113 & 41 \\ \hline
		\end{tabular}
	\end{minipage}
\end{table}

\section{Conclusion}\label{sec:con_fut}
We discuss application of preconditioned iterative methods for solving the sequences of large sparse linear systems in a class of parametric model order reduction algorithms (for reducing second-order linear dynamical systems). We specifically focus on Robust parametric model order reduction algorithm of \cite{lihong2014}, however, show that our techniques are applicable to other such algorithms as well. 

The choice of the iterative method is problem dependent, however, because the multiple right-hand-sides of the linear systems are available here together, we propose the use of block variant of the chosen iterative method. We show use of the SPAI preconditioner because it is inherently parallel, and especially useful in solving exponentially increasing linear systems sizes that arise here. 

The linear systems here change slightly during the model order reduction process. Exploiting this, we propose a technique to cheaply update the SPAI preconditioner. The novelty here is that our update exploits the slowly changing behavior of the model parameters, which was not done earlier. For a model problem, while using GCRO as the underlying iterative solver, we show that using a block GCRO saves 80\% of computation time over its non-block version. Further, by using SPAI update, we  save 70\% of the time over simple SPAI preconditioner.  

In future, we plan to apply preconditioned iterative methods to other model order reduction algorithms (for those higher than second order as well as bilinear and non-linear dynamical systems). We also plan to look at the possibility of developing general preconditioners for each of these types of model order reduction algorithms.

\section*{References}

\bibliography{references}
}
\section*{Appendix I} 
\myblue{In \cite{Ionita2014}, PMOR-L is tested on the example from \cite{HODEL1996205,Penzl1999}. This  corresponds to a second order linear dynamical system, which is equivalent to first order linear dynamical systems given in (\ref{eq1:fst-order-ds}). We now look at arising linear system matrices as follows:
\begin{align*}
\mathbb{E}(p_j) = I, \ \mathbb{A}(p_j) = diag\left( A_1(p_j), A_2, A_3, A_4 \right) \ \text{for} \ j = 1, \ldots w,
\end{align*}
where $A_1(p_j) = \begin{bmatrix}
-1 & p_j \\
-p_j & -1
\end{bmatrix}$, $A_2 = \begin{bmatrix}
-1 & 200 \\
-200 & -1
\end{bmatrix}$, $A_3 = \begin{bmatrix}
-1 & 400 \\
-400 & -1
\end{bmatrix}$, and $A_4 = diag \left(1, 2, \ldots 1000\right)$. Thus,
\begin{align*}
A(s, p_1) 
       = & sI -   \begin{bmatrix}
        A_1(p_1) & & & \\
        & A_2 & & \\
        & & A_3 & \\
        & & & A_4
        \end{bmatrix}\\ 
        & \text{and} \\
        A(s, p_j) 
        = & sI -   \begin{bmatrix}
        A_1(p_j) & & & \\
        & A_2 & & \\
        & & A_3 & \\
        & & & A_4
        \end{bmatrix}.
\end{align*}
Now, expressing $A(s,p_j)$ in term of $A(s,p_1)$ we get
\begin{align*}
A(s,p_j) = & A(s,p_1) +   \begin{bmatrix}
A_1(p_1) & & & \\
& A_2 & & \\
& & A_3 & \\
& & & A_4
\end{bmatrix} - \begin{bmatrix}
A_1(p_j) & & & \\
& A_2 & & \\
& & A_3 & \\
& & & A_4
\end{bmatrix} \\
= & A(s,p_1) + diag \left((A_1(p_1)-A_1(p_j)), 0, \ldots, 0 \right).
\end{align*}
Thus, preconditioner update corresponding to parameter values $p_j$ can be easily applied.}

\section*{Appendix II}
\myblue{In \cite{sonN2017}, PBTMR is tested on a heat model. Next, we look at the arising linear system matrices for this model.
	\begin{align*}
	A(p_j) = - \mathbb{E}(p_j) \otimes \mathbb{A}(p_j) -\mathbb{A}(p_j) \otimes \mathbb{E}(p_j), \ \text{for}\ j = 1, \ldots, w, 
	\end{align*}
	with $\mathbb{E}(p_j) = E, \ \mathbb{A}(p_j) = \tilde{p}_1^j A_1 + \tilde{p}_2^j A_2 + \tilde{p}_3^j A_3 + \tilde{p}_4^j A_4 + A_5$, and $p_j = [\tilde{p}_1^j, \tilde{p}_2^j, \tilde{p}_3^j, \tilde{p}_4^j]^{T}$. Thus, 
	}
	\myblue{ 
		\begin{align*}\label{eq:lap_ls_1}
		\begin{split}
		A(p_1) 
		& = -\left(E \otimes  (\tilde{p}_1^1 A_1+ \tilde{p}_2^1 A_2 + \tilde{p}_3^1 A_3 + \tilde{p}_4^1 A_4) + E \otimes A_5  + (\tilde{p}_1^1 A_1+ \tilde{p}_2^1 A_2 + \tilde{p}_3^1 A_3 + \tilde{p}_4^1 A_4) \otimes E + A_5 \otimes E \right), \\
	    & \text{and} \\
		A(p_j) 
		& = -\left(E \otimes  (\tilde{p}_1^j A_1+ \tilde{p}_2^j A_2 + \tilde{p}_3^j A_3 + \tilde{p}_4^j A_4) + E \otimes A_5  + (\tilde{p}_1^j A_1+ \tilde{p}_2^j A_2 + \tilde{p}_3^j A_3 + \tilde{p}_4^j A_4) \otimes E + A_5 \otimes E \right).
		\end{split}	 
		\end{align*} 
		Now, expressing $A(p_j)$ in term of $A(p_1)$ we get
		\begin{align*}
		A(p_j) = & \ A(p_1) + E \otimes  (\tilde{p}_1^1 A_1+ \tilde{p}_2^1 A_2 + \tilde{p}_3^1 A_3 + \tilde{p}_4^1 A_4) + (\tilde{p}_1^1 A_1+ \tilde{p}_2^1 A_2 + \tilde{p}_3^1 A_3 + \tilde{p}_4^1 A_4) \otimes E \\ & -  E \otimes  (\tilde{p}_1^j A_1+ \tilde{p}_2^j A_2 + \tilde{p}_3^j A_3 + \tilde{p}_4^j A_4) - (\tilde{p}_1^j A_1+ \tilde{p}_2^j A_2 + \tilde{p}_3^j A_3 + \tilde{p}_4^j A_4) \otimes E,\\
		& = A(p_1) + (EA)_j,
		\end{align*}
		where $(EA)_j = E \otimes  (\tilde{p}_1^1 A_1+ \tilde{p}_2^1 A_2 + \tilde{p}_3^1 A_3 + \tilde{p}_4^1 A_4) + (\tilde{p}_1^1 A_1+ \tilde{p}_2^1 A_2 + \tilde{p}_3^1 A_3 + \tilde{p}_4^1 A_4) \otimes E -  E \otimes  (\tilde{p}_1^j A_1+ \tilde{p}_2^j A_2 + \tilde{p}_3^j A_3 + \tilde{p}_4^j A_4) - (\tilde{p}_1^j A_1 +  \tilde{p}_2^j A_2 + \tilde{p}_3^j A_3 + \tilde{p}_4^j A_4) \otimes E$.
		Thus preconditioner update can be easily applied.}

\end{document}